\documentclass[psamsfonts]{amsart}

\newtheorem{theorem}{Theorem}[section]
\newtheorem{lemma}[theorem]{Lemma}
\newtheorem{proposition}[theorem]{Proposition}
\newtheorem{corollary}[theorem]{Corollary}

\theoremstyle{definition}

\newtheorem{example}[theorem]{Example}
\newtheorem{question}[theorem]{Question}

\newtheorem{remark}[theorem]{Remark}

\theoremstyle{remark}

\numberwithin{equation}{section}

\newcommand{\NN}{\mathbb{N}}
\newcommand{\ZZ}{\mathbb{Z}}

\newcommand{\CC}{\mathbb{C}}

\newcommand{\PP}{\mathbb{P}}
\renewcommand{\AA}{\mathbb{A}}

\newcommand  {\shF}     {\mathcal{F}}
\newcommand  {\shG}     {\mathcal{G}}

\newcommand  {\shL}     {\mathcal{L}}
\newcommand  {\shR}     {\mathcal{R}}
\newcommand  {\shS}     {\mathcal{S}}

\newcommand  {\foa}     {\mathfrak{a}}

\newcommand  {\Det}    {\operatorname{Det}}

\newcommand  {\dual}    {\vee}

\newcommand  {\lra}     {\longrightarrow}

\renewcommand{\O}       {\mathcal{O}}

\newcommand  {\Proj}    {\operatorname{Proj}}

\newcommand  {\ra}      {\rightarrow}

\newcommand  {\rank}    {\operatorname{rank}}

\newcommand  {\Spec}    {\operatorname{Spec}}

\newcommand {\pasoclo} {\star}

\def\mydate{\number\day\space\ifcase\month \or January\or February\or March\or April\or May\or
June\or July\or August\or September\or October\or November\or
December\fi \space\number\year}

\usepackage{amscd}
\usepackage{amssymb}

\begin{document}

\title[Tight closure and elliptic curves]
{Tight closure and plus closure for cones over elliptic curves}

% Remove or comment out any unused author tags.
% author one information

\author[Holger Brenner]{Holger Brenner}
\address{Department of Pure Mathematics, University of Sheffield,
Hicks Building, Hounsfield Road, Sheffield S3 7RH, United Kingdom}
%\curraddr{}
\email{H.Brenner@sheffield.ac.uk}

%\thanks{}

%\date{}

% at present the "communicated by" line appears only in ERA and PROC
%\commby{}

\dedicatory{\mydate}

\begin{abstract}
We characterize the tight closure of a homogeneous primary ideal
in a normal homogeneous coordinate ring over an elliptic curve
by a numerical condition
and we show that it is in positive characteristic
the same as the plus closure.
\end{abstract}

\maketitle


%===========================================================
\section*{Introduction}

Let $\foa =(f_1 , \ldots ,f_n) \subseteq R$ denote an ideal in a Noetherian domain
over a field of positive characteristic $p$.
Hochster and Huneke introduced the notion of the tight closure of the ideal $\foa$,
which is given by
$$\foa^* = \{f \in R:\, \exists c \neq 0 \mbox{ such that }
cf^q \in (f_1^q, \ldots ,f_n^q) \mbox{ for all powers } q=p^{e}\}$$

One of the basic open questions in tight closure theory
is the problem whether the tight closure of an ideal $\foa$
in a domain $R$ of positive characteristic is just the contraction
$R \cap \foa R^+$ from the absolute integral closure $R^+$ of $R$.
Hochster calls this a tantalizing question. A positive answer would
imply that tight closure commutes with localization.
The best result so far is given by
the theorem of Smith \cite{smithparameter}, \cite[Theorem 7.1]{hunekeparameter}
which states that for parameter ideals
the tight closure and the plus closure are the same.

The general question is open even in the case of a two-dimensional
normal graded domain (in a regular ring every ideal is tightly closed,
so there is no problem).
The domain $R=K[x,y,z]/(F)$ for the
Fermat polynomial $F=x^3+y^3+z^3$
is a standard example in tight closure theory
and has been intensively studied (\cite{mcdermott},\cite{singh}),
but even in this simple looking
example neither the answer to the question is known
nor is it clear how to compute the tight closure of a given ideal.

In this paper we study the normal homogeneous coordinate ring $R$ of an
elliptic curve $Y$
over an algebraically closed field $K$.
That is $R$ is a normal standard-graded two-dimensional $K$-algebra
such that $Y=\Proj \, R$ is an elliptic curve.
This contains in particular the case $R=K[x,y,z]/(F)$,
where $F$ is a homogeneous polynomial of degree $3$ such that $\Spec\, R$
is non-singular outside the origin.
Our main result is that for an $R_+$-primary homogeneous ideal $\foa \subseteq R$
the tight closure and the plus closure are the same.
This follows from a numerical criterion which holds for both
closure operations.

We obtain these results by applying the geometric method
which we developed in \cite{brennertightproj}. For tight closure data
consisting of homogeneous generators $f_1, \ldots ,f_n$
of an $R_+$-primary ideal in a graded ring $R$ and another homogeneous element
$f_0$ we construct a projective bundle together with a projective
subbundle of codimension one over $Y=\Proj\, R$.
The questions whether $f_0 \in (f_1, \ldots , f_n)^*$ and
$f_0 \in (f_1, \ldots , f_n)^+$
translates then to questions on the complement of the subbundle:
whether it is non-affine and whether it contains projective curves
(for $\dim\, R=2$).
Therefore we can work in a projective geometric setting.
We recall the construction and the necessary facts in section \ref{graded}.

Since we are concerned with homogeneous coordinate rings over
elliptic curves, we are in a very favorable situation: The vector
bundles on an elliptic curve have been completely classified by
Atiyah \cite{atiyahelliptic}, leading to further results on the
ampleness and on the behavior of cohomology classes, which we
recollect and extend in section \ref{elliptic} for our needs.

In section \ref{ellipticaffin}
we give for an extension
$0 \ra \shR \ra \shR' \ra \O_Y \ra 0$ of locally free sheaves a numerical criterion
for the affineness of $\PP(\shR^{' \dual}) - \PP(\shR^\dual)$ in terms
of the cohomology class $\in H^1(Y, \shR)$ of the extension
and the degrees of the indecomposable components
of $\shR$ (Theorem \ref{numkritaffin}).
In positive characteristic we obtain the same numerical condition
for the property that every projective curve meets the projective subbundle
$\PP(\shR^\dual)$
(Theorem \ref{numkritcurve}).

In section \ref{tightplussection} we derive from these geometric results the corresponding
statements for tight closure and plus closure for
primary homogeneous ideals in a homogeneous coordinate ring over an elliptic
curve.

%===========================================================

\section{Projective bundles corresponding to tight closure problems}
\label{graded}

In this section we recall how graded tight closure problems in a graded ring
$R$ translate to
problems on projective bundles and subbundles over $\Proj\, R$.
Let's fix our notations.
By a vector bundle $V$ on a scheme $Y$ we mean a geometric vector bundle.
We denote the locally free sheaf of sections in $V$ by $\shR$
and its dual, the sheaf of linear forms, by $\shF$.
Hence $V= \Spec\, \oplus_{n \geq 0}\, S^n(\shF)$,
and $\PP(V) = \PP(\shF) = \Proj \, \oplus_{n \geq 0}\, S^n(\shF)$.
For a number $k$ we denote
by $\AA_Y(k)$ the (geometric) line bundle
$\Spec \, \oplus_{n \geq 0} \O_Y(nk)$
with sheaf of sections $\O_Y(-k)$.
We consider sometimes the geometric realization
$\AA_Y(k)= D_+(R_+) \subset \Proj\, R[T] $, where
$ \deg \, T =  -k$.

Let $K$ denote an algebraically closed field and let $R$ be a
standard $\NN$-graded $K$-algebra, that is $R_0=K$ and $R$ is
generated by finitely many elements of degree one. Set $Y= \Proj \,
R$. Let $f_i$ be homogeneous primary elements of $R$ of degrees
$d_i$, that is the $D_+(f_i)$ cover $Y$. Fix a number $m \in \ZZ$.
Let $A=R[T_1, \ldots , T_n]/(f_1T_1 + \ldots +f_nT_n)$ be graded by
$\deg \, T_i = m-d_i$ (maybe negative). Then the open subset $\Spec
\, A \supset D(R_+A)$ is a vector bundle over $D(R_+) \subset
\Spec\, R$.

\begin{proposition}
\label{buendelalsproj}
Let $R$ be a standard-graded $K$-algebra and let
$f_1,\ldots,f_n$ be homogeneous primary elements.
Let $d_i = \deg \, f_i $ and let $m \in \ZZ$.
Then the following hold.

\renewcommand{\labelenumi}{(\roman{enumi})}
\begin{enumerate}

\item
$$ \Proj\, R[T_1,\ldots,T_n]/(\sum_{i=1}^n f_iT_i) \, \,
\supset \, \, D_+(R_+)
\longrightarrow \Proj\, R$$
is a vector bundle of rank $n-1$ over $Y=\Proj\, R$, which we denote by $V(-m)$ .

\item
There exists an exact sequence of vector bundles
$$0 \longrightarrow V(-m) \longrightarrow
\AA_Y(d_1 -m) \times_Y \ldots \times_Y \AA_Y(d_n - m)
\stackrel{\sum f_i}{\lra} \AA_Y(-m) \longrightarrow 0 \, .$$

\item
We have $\Det \, V(-m) \, \cong \,  \AA_Y(\sum_{i=1}^n d_i -(n-1)m) $, and
$\deg \, V( - m) =(\sum_{i=1}^n d_i -(n-1)m) \deg \, \AA_Y(1)$.

\item
We have $V (-m')= V( - m) \otimes \AA_Y(m-m')$,
and $\PP(V( - m))$ does not depend on the degree $m$.

\end{enumerate}
\end{proposition}

\proof
See \cite[Proposition 3.1]{brennertightproj}.
\qed

\begin{proposition}
\label{forcingsequence2}
Let $R$ be a standard-graded $K$-algebra, let
$f_1,\ldots,f_n$ be homogeneous primary elements and let
$f_0 \in R$ be also homogeneous.
Let $d_i = \deg \, f_i $, $m \in \ZZ$ and set $ \deg \, T_i =m -d_i$.
Let
$$V( - m)=D_+(R_+) \subset \Proj\, R[T_1,\ldots,T_n]/(\sum_{i=1}^n f_iT_i)
\, \, \mbox{ and }$$
$$V'( - m)= D_+(R_+) \subset \Proj\, R[T_0,\ldots,T_n]/(\sum_{i=0}^n f_iT_i)$$
be the vector bundles on $Y=\Proj\, R$ due to {\rm \ref{buendelalsproj}}
Then the following hold.

\renewcommand{\labelenumi}{(\roman{enumi})}
\begin{enumerate}

\item
There is an exact sequence of vector bundles on $Y$,
$$0 \longrightarrow V( - m) \longrightarrow V'( - m)
\stackrel{T_0}{\longrightarrow} \AA_Y(d_0-m) \longrightarrow 0 \, .$$

\item
The embedding $\PP(V( - m)) \hookrightarrow \PP(V'( - m))$ does not depend
on $m$.

\item
Let $E$ be the hyperplane section on $\PP(V')$
corresponding to the relative very ample
invertible sheaf $\O_{\PP(V')} (1)$ {\rm (}depending of the degree{\rm )}.
Then we have the linear equivalence of divisors
$\PP(V) \sim E +(m-d_0) \pi^*H$,
where $H$ is the hyperplane section of $\, Y$.
If $m=d_0$, then $\PP(V)$ is a hyperplane section.

\end{enumerate}
\end{proposition}
\proof
See \cite[Proposition 3.4]{brennertightproj}.
\qed

\begin{remark}
\label{forcingdef}
We call the sequence in \ref{forcingsequence2} (i)
the forcing sequence. We often skip the number $m$
(the total degree)
and denote the situation $\PP(V) \hookrightarrow \PP(V')$
by $\PP(f_1, \ldots, f_n; f_0)$.
$\PP(V)=\PP(f_1, \ldots, f_n)$ is
called the forcing divisor or forcing subbundle.
The complement $\PP(V') -\PP(V)$ plays a crucial role in our method,
since it is isomorphic to the $\Proj$
of the so called forcing algebra
$R[T_1, \ldots, T_n]/(f_1T_1 + \ldots +f_nT_n +f_0)$
(suitable graded).

We denote the sheaves of sections in $V( - m)$ (and $V'( - m)$) by
$\shR (m)$ (and $\shR'(m)$). This is the sheaf of syzygies (or
relations) for the elements $f_1, \ldots, f_n$ of total degree $m$.
We denote the corresponding sheaves of linear forms by $\shF(-m)
=\shR (m)^\dual$ and $\shF'(-m)$. Note that for $m=d_0$ the forcing
sequence $0 \ra \shR(m) \ra \shR'(m) \ra \O_Y \ra 0$ corresponds to
a cohomology element $c \in H^1(Y,\shR(m))$.
\end{remark}

\medskip
The containment of a homogeneous element in the ideal,
in the tight closure and in the plus closure of the ideal
is expressed in terms of the projective bundles in the following way.
In the case of characteristic zero,
the notion of plus closure does not make much sense
and we work rather with solid closure than with tight closure,
see \cite{hochstersolid}.

\begin{lemma}
\label{trivialtwo}
In the situation of {\rm \ref{forcingsequence2}} the following
are equivalent.

\renewcommand{\labelenumi}{(\roman{enumi})}
\begin{enumerate}

\item
$f_0 \in (f_1, \ldots ,f_n)$.

\item
There is a section $Y \ra \PP(V')$ disjoined to $\PP(V) \subset \PP(V')$.

\item
The forcing sequence
$0 \ra V(-m) \ra V'(-m) \ra \AA_Y(d_0-m) \ra 0 $ splits.

\item
Suppose $m=d_0$.
The corresponding cohomology class in $H^1(Y,\shR (m))$ vanishes.
\end{enumerate}
\end{lemma}

\proof
See \cite[Lemma 3.7]{brennertightproj}.
\qed

\begin{proposition}
\label{geokrittightplus} Let $R$ be a normal standard-graded
$K$-algebra of dimension $2$, let $f_1, \ldots ,f_n \in R$ be
primary homogeneous elements and let $f_0$ be another homogeneous
element. Let $V$ and $V'$ be as in {\rm \ref{forcingsequence2}}.
Then $f_0 \in (f_1, \ldots ,f_n)^*$ if and only if $\PP(V') -
\PP(V)$ is not affine.

Furthermore, if the characteristic of $K$ is positive,
the following are equivalent.
\renewcommand{\labelenumi}{(\roman{enumi})}
\begin{enumerate}

\item
$f_0 \in (f_1, \ldots ,f_n)^{+{\rm gr}}$, i.e. there exists
a finite graded extension $R \subseteq R'$ such that
$f_0 \in (f_1, \ldots ,f_n)R'$.

\item
There exists a smooth projective curve $X$
and a finite surjective morphism $g: X \ra Y$ such that
the pull back $g^*\PP(V')$
has a section not meeting $g^*\PP(V)$.

\item
There exists a curve $X \subset \PP(V')$
which does not intersect $\PP(V)$.

\end{enumerate}
\end{proposition}
\proof
See \cite[Lemmata 3.9 and 3.10]{brennertightproj}.
\qed

\section{Vector bundles over elliptic curves}
\label{elliptic}

We gather together some results on vector bundles over elliptic curves.
Recall that a locally free sheaf $\shF$ on a scheme $Y$ is called ample
if the invertible sheaf $\O_{\PP(\shF)}(1)$ on the projective bundle
$\PP(\shF)$ is ample.
If $\shF'$ is ample and $0 \ra \O_Y \ra \shF' \ra \shF \ra 0$
is a short exact sequence, then $\PP(\shF) \subset \PP(\shF')$
is an ample divisor, hence its complement is affine.
The following theorem of Gieseker-Hartshorne
gives a numerical criterion for ample bundles over an elliptic curve.

\begin{theorem}
\label{amplekrit}
Let $Y$ denote an elliptic curve over an algebraically closed field $K$
and let $\shF$ denote a locally free sheaf.
Then $\shF$ is ample if and only if the degree of every
indecomposable summand of $\shF$ is positive.
\end{theorem}

\proof
See \cite[Prop. 1.2 and Theorem 1.3]{hartshorneamplecurve} or
\cite[Theorem 2.3]{giesekerample}.
\qed

\medskip
The following theorem is a generalization of a theorem of Oda.

\begin{theorem}
\label{odageneral}
Let $Y$ denote an elliptic curve over an algebraically closed field $K$
and let $\shR$ denote an indecomposable locally free sheaf of negative degree.
Let $g: X \ra Y$ be a finite dominant map, where $X$ is another smooth
projective curve.
Then $H^1(Y, \shR) \lra H^1(X, g^* \shR)$ is injective.
\end{theorem}

\proof
Let $0 \neq c \in H^1(Y, \shR)$ be a non zero class and consider
the corresponding extension
$$ 0 \lra \shR \lra \shR' \lra \O_Y \lra 0 \, .$$
Let $\shF$ and $\shF'$ denote the dual sheaves and let
$\PP(\shF) \hookrightarrow \PP(\shF')$ be the corresponding
projective subbundle.
The indecomposable sheaf $\shF$ is of positive degree,
hence ample due to \ref{amplekrit}.
Then also $\shF'$ is ample due to
\cite[Proposition 2.2]{giesekerample} and \ref{amplekrit}.
Hence $\PP(\shF)$ is an ample divisor on $\PP(\shF')$ and
its complement is affine. This property is preserved under the finite
mapping $X \ra Y$, therefore the complement of
$\PP(g^*(\shF)) \subset \PP(g^*(\shF'))$ is affine and there cannot
be projective curves in the complement. Hence the pulled back sequence
does not split and $g^*(c) \neq 0$.
\qed

\begin{remark}
Oda proved this statement only for the Frobenius morphism in
positive characteristic, see \cite[Theorem 2.17]{odaelliptic},
and Hartshorne used this theorem to prove the numerical criterion for
ampleness. The proof of this criterion by Gieseker
in \cite{giesekerample}
however is independent
of the theorem of Oda.
\end{remark}

\medskip
A kind of reverse to \ref{odageneral} is given by the following lemma.

\begin{lemma}
\label{pvanish}
Let $K$ be an algebraically closed
field of positive characteristic $p$ and let $Y$ be an elliptic
curve.
Let $\shG$ be an indecomposable locally free sheaf on $Y$ of degree $\geq 0$.
and let $c \in H^1(Y, \shG)$ be a cohomology class.
Then there exists a finite curve
$g: X \ra Y $ such that $g^*(c) \in H^1(X, g^*(\shG))$ is zero.
\end{lemma}

\proof
In fact we will show that the multiplication mappings $[p^{e}]:Y \ra Y$
have the stated property.
If $\deg \, \shG > 0$, then
$H^1(Y, \shG) \cong H^0 (Y, \shG^\dual)=0$ due to
\cite[Lemma 1.2]{hartshorneamplecurve} and there is nothing to prove.
The same is true for an invertible sheaf $\neq \O_Y$ of degree $0$.
For $\shG =\O_Y$, the multiplication map
$[p]: Y \ra Y$ induces the zero map on
$H^1(Y, \O_Y)$, this follows for example
from \cite[Corollary 5.3]{silvermanelliptic} and
\cite[\S 13 Corollary 3]{mumfordabelian}.

Now we do induction on the rank, and
suppose that $r=\rank\, \shG \geq 2$ and $\deg \, \shG =0$.
Due to the classification of Atiyah (see \cite[Theorem 5]{atiyahelliptic})
we may write
$\shG = \shF_r \otimes \shL$, where $\shL$ is an invertible sheaf
of degree $0$ and where
$\shF_r$ is the unique sheaf of rank $r$ and degree zero with
$\Gamma(Y, \shF_r) \neq 0$.
In fact, for these sheaves we know that $H^0(Y, \shF_r)$ and $H^1(Y, \shF_r)$
are one-dimensional and there exists a non-splitting
short exact sequence
$$ 0 \lra \O_Y \lra \shF_{r} \lra \shF_{r-1} \lra  0 \, .$$
This gives the sequence
$ 0 \ra \shL  \ra \shG \ra \shF_{r-1} \otimes \shL \ra 0$.
Let $c \in H^1(Y, \shG)$. Then the image of this class in
$H^1(Y, \shF_{r-1} \otimes \shL)$ is zero after applying
$[p^{e}]$ and comes then from an element in $H^1(Y, \shL)$, which itself is
zero after applying $[p]$ once more.
\qed

\section{A numerical criterion for subbundles to have affine complement}
\label{ellipticaffin}

In this section we investigate subbundles $\PP(\shF) \subset \PP(\shF')$
of codimension one over an elliptic curve with respect to the properties
which are interesting from the tight closure and plus closure point of view:
Is the complement affine? Does it contain projective curves?

\begin{lemma}
\label{affinlemma}
Let $Y$ denote a scheme
and let
$\shR$ and $\shS$
be locally free sheaves on $Y$ and let
$\varphi: \shR \ra \shS$ be a morphism.
Let $c \in H^1(Y,\shR)$ with corresponding extension
$0 \ra \shR \ra \shR' \ra \O_Y \ra 0$
and let $d \in H^1(Y, \shS)$ be its image
with corresponding extension
$0 \ra \shS \ra \shS' \ra \O_Y \ra 0$.
Let $V$ and $W$ denote the corresponding vector bundles
{\rm (}$V=\Spec \, S(\shR^\dual)$ etc.{\rm )}
If then $\PP(W') - \PP(W)$ is affine, then also
$\PP(V')-\PP(V)$ is affine.

\end{lemma}
\proof
First note that we have a mapping
$\varphi': V' \ra W'$ compatible with the
extensions and with $\varphi$,
see \cite[Ch. 3 Lemma 1.4]{maclanehomology}.
The induced rational mapping
$\PP(V') \ra \PP(W')$ is defined outside the kernel of $\varphi '$.
Locally these mappings on the vector bundles look like
$$
\begin{CD}
\AA^r    @>  >> \AA^r  \times \AA   \\
 \varphi  @VVV     @VVV  \varphi \times id \\
\AA^s    @>  >> \AA^s  \times \AA  \, \, .
\end{CD}
$$
The line on $\PP(V')$ corresponding to the point $v=(v,t)$
does not lie on the subbundle $\PP(V)$ for $t \neq 0$. For these points
the rational mapping is defined and the image point $(\varphi(v),t)$
does not lie on the subbundle $\PP(W)$.
Hence we have an affine morphism
$$\PP(V')-\PP(V) \, \lra \,  \PP(W') - \PP(W) $$
and
if $\PP(W') - \PP(W)$ is affine, then
$\PP(V')-\PP(V)$ is affine as well.
\qed

\medskip
Our first main result is the following numerical characterization
for the complement of a projective subbundle to be affine.

\begin{theorem}
\label{numkritaffin}
Let $K$ be an algebraically closed field and let
$Y$ denote an elliptic curve.
Let $\shR$ be a locally free sheaf on $Y$ of
rank $r$ and let $\shR = \shR_1 \oplus \ldots \oplus \shR_s$
be the decomposition into indecomposable locally free sheaves.
Let $c \in H^1(Y, \shR)$ and let
$0 \ra \shR \ra \shR' \ra \O_Y \ra 0$ be the corresponding
extension and let $\PP(\shF) \subset \PP(\shF')$ be the corresponding
projective bundles. Then the following are equivalent.

\renewcommand{\labelenumi}{(\roman{enumi})}
\begin{enumerate}

\item
There exists
$1 \leq j \leq s$ such that $\deg \shR_j < 0$
and $c_j \neq 0$, where $c_j $ denotes the component of $c$
in $H^1(Y, \shR_j)$.

\item
The complement $\PP(\shF') - \PP(\shF)$ is affine.
\end{enumerate}

\end{theorem}
\proof
(i) $\Rightarrow$ (ii).
Suppose that $j$ fulfills the statement in the numerical criterion.
Consider the projection $p_j :  \shR  \ra \shR_j$,
where the corresponding cohomological map sends $c$ to $c_j$.
Now $\shF_j= \shR_j^\dual$ is indecomposable of positive degree, hence ample
due to \ref{amplekrit}. Since $\shF_j'$ is a non-trivial
extension of $\shF_j$, it is also ample on the elliptic curve $Y$,
see \cite[Proposition 2.2 and Theorem 2.3]{giesekerample} (this is true
for every curve in characteristic zero, see
\cite[Theorem 2.2]{giesekerample}).
But then the complement of the (hypersection) divisor
$\PP(\shF_j) \subset \PP(\shF_j')$
is affine, hence $\PP(\shF')-\PP(\shF)$
is affine due to \ref{affinlemma}.

(ii) $ \Rightarrow$ (i).
Suppose to the contrary that
$\deg \, \shR_j \geq 0$ or $c_j=0$ holds for every $j$.
Since $H^1(Y, \shR_j)=H^0(Y, \shR_j^\dual)=0$ for $\deg \, \shR_j >0$
due to \cite[Lemma 1.1]{hartshorneamplecurve}, we may assume that
$\deg \, \shR_j = 0$ or $c_j=0$ holds for every $j$.
Due to \ref{affinlemma} we may forget the components with $c_j=0$
and hence assume that $\deg\, \shR_j=0$ for every component.

We claim that $\PP(\shF)$ is a numerically effective divisor
on $\PP(\shF')$.
Since indecomposable sheaves on elliptic curves are semistable and
since all the components of $\shF$ have degree $0$ it follows that
$\shF$ is semistable.
Therefore $\deg \, \shL \geq 0$
for every quotient invertible sheaf $\shL$ of $\shF$.

A semistable sheaf on an elliptic curve stays semistable
after applying a finite dominant morphism $X \ra Y$,
see \cite[Proposition 5.1]{miyaokachern},
hence also the degree of a quotient invertible sheaf is non-negative
on every curve $X \ra Y$.
It then follows that also the degree of an invertible quotient
sheaf of $\shF'$ is nonnegative on every curve,
and this means that the intersection of $\PP(\shF)$
with any curve is nonnegative.

Since the degree of our numerically effective divisor is zero,
it follows by the Kodaira Lemma \cite[Lemma 2.5.7]{beltramettiadjunction}
that it is not big.
Therefore its complement cannot be affine.
\qed

\medskip
We shall show now that the same numerical criterion
holds in positive characteristic for the (non-)existence of projective curves
inside $\PP(\shF') -\PP(\shF)$.

\begin{theorem}
\label{numkritcurve}
Let $K$ be an algebraically closed field of positive characteristic and let
$Y$ denote an elliptic curve.
Let $\shR$ be a locally free sheaf on $Y$ of
rank $r$ and let $\shR = \shR_1 \oplus \ldots \oplus \shR_s$
be the decomposition in indecomposable locally free sheaves.
Let $c \in H^1(Y, \shR)$ and let
$0 \ra \shR \ra \shR' \ra \O_Y \ra 0$ be the corresponding
extension and let $\PP(\shF) \subset \PP(\shF')$ be the corresponding
projective bundles.
Then the following are equivalent.

\renewcommand{\labelenumi}{(\roman{enumi})}
\begin{enumerate}

\item
There exists
$1 \leq j \leq s$ such that $\deg \shR_j < 0$
and $c_j \neq 0$, where $c_j $ denotes the component of $c$
in $H^1(Y, \shR_j)$.

\item
The sequence
$0 \ra \shR \ra \shR' \ra \O_Y \ra 0$ does not split after
a finite dominant morphism $X \ra Y$, where $X$ is another projective
curve.

\item
The subbundle $\PP(\shF) \subset  \PP(\shF')$ intersects every curve
in $\PP(\shF')$ positively.

\end{enumerate}

\end{theorem}
\proof

(i) $\Rightarrow$ (ii).
Suppose that the sequence splits under the finite morphism
$g:X \ra Y$. Then $g^*(c_j)=0$ on $X$ and from
\ref{odageneral}
we see that $\deg \, \shR_j \geq 0$ or $c_j=0$

(ii) $\Rightarrow$ (iii).
If there exists a curve $C$ on
$\PP(\shF')$ not meeting $\PP(\shF)$,
then it dominates the base. Let $X$ be the normalization of $C$ and
let $g: X \ra Y$ the finite dominant mapping. Then
$g^*\PP(\shF') \ra X$ has a section not meeting $g^*\PP(\shF)$
and then the sequence splits on $X$.

(iii) $\Rightarrow $ (i).
Suppose to the contrary that for all the indecomposable components
of $\shR$ with negative degree $c_j =0$ holds.
For every component $\shR_j$ with $\deg \, R_j \geq 0$ there exists
due to \ref{pvanish}
a finite curve $g_j: X_j \ra Y$ such that
$g_j^*(c_j) \in H^1(X_j, g_j^* \shR_j)$
is zero.
Putting these curves together we find a curve $g: X \ra Y$
such that $g^* (c)=0$.
Thus the sequence splits on $X$ and this gives a projective curve
in $\PP(\shF') - \PP(\shF)$.
\qed

\begin{corollary}
\label{geokritaffin}
Let $K$ be an algebraically closed field of positive characteristic and let
$Y$ denote an elliptic curve.
Let $\shR$ be a locally free sheaf on $Y$,
let $c \in H^1(Y, \shR)$ and let
$0 \ra \shR \ra \shR' \ra \O_Y \ra 0$ be the corresponding
extension.
Let $\PP(\shF) \subset \PP(\shF')$ be the corresponding
projective bundles.
Then $\PP(\shF') - \PP(\shF)$ is affine if and only if
it contains no
projective curve.
\end{corollary}
\proof
This follows from \ref{numkritaffin} and \ref{numkritcurve}, since for
both properties the same numerical criterion holds.
\qed

\begin{remark}
The corollary does not hold in characteristic zero.
The sequence $0 \ra \O_Y \ra \shF_2 \ra \O_Y \ra 0$ yields a section
in a ruled surface whose complement is not affine, but Stein (over $\CC$),
hence it does not contain projective curves.
\end{remark}

\begin{question}
We will derive from the corollary in the next section that the tight closure
of a primary homogeneous ideal
is the same as its plus closure in the normal homogeneous coordinate ring
over an elliptic curve. A natural generalization of the question
whether $\foa^* = \foa^+$
is the following question.

Let $0 \ra \shR \ra \shR' \ra \O_Y \ra 0$ be an exact sequence
of locally free sheaves
on a (smooth) projective variety of dimension $d$ over an
(algebraically closed) field of
positive characteristic (!).
Let $\PP(\shF) \subset \PP(\shF')$
be the corresponding bundles
and suppose that
$\PP(\shF)$ meets every subvariety of dimension $d$.
Is then the cohomological dimension
$cd (\PP(\shF') - \PP(\shF)) <d$?
If the sequence is a forcing sequence and if $\rank\, \shR =d$,
then this is true due to the parameter theorem of Smith.
\end{question}

\medskip
Tight closure theory yields a lot of further questions concerning
vector bundles on projective varieties. Let me just mention the
following.

\begin{question}
Let
$0 \ra \shR \ra \shR' \ra \O_Y \ra 0$ be an exact sequence
of locally free sheaves on $\PP^d$ (any characteristic),
$\PP(\shF) \subset \PP(\shF')$.
Suppose that the mapping
$H^d(\PP^d, \O_{\PP^d}(k)) \ra
H^d (\PP(\shF')- \PP(\shF), \pi^*\O_{\PP^d}(k))$
is injective for $k \leq 0$ (then of course for all $k$).

Does the sequence split?
This is true for forcing sequences.
\end{question}

\section{Numerical criteria for tight closure and plus closure}
\label{tightplussection}

We are now in the position to draw the consequences to tight closure
(solid closure in characteristic $0$)
and plus closure in a normal homogeneous coordinate ring of an elliptic curve,
that is $R$ is a two-dimensional normal standard-graded domain such that
$\Proj\, R$ is an elliptic curve.

\begin{corollary}
\label{numkrittight} Let $K$ be an algebraically closed field and
let $R$ be a normal homogeneous coordinate ring over the elliptic
curve $Y= \Proj \, R$. Let $f_1, \ldots ,f_n$ be homogeneous
generators of an $R_+$-primary homogeneous ideal in $R$ and let $m
\in \NN$ be a number. Let $\shR(m)$ be the corresponding locally
free sheaf of syzygies of total degree $m$ on $Y$ of rank $n-1$ and
let $\shR(m) = \shR_1 \oplus \ldots \oplus \shR_s$ be the
decomposition into indecomposable locally free sheaves. Let $f_0 \in
R$ be another homogeneous element of degree $m$ and let $c \in
H^1(Y, \shR(m))$ be the corresponding class. Then the following are
equivalent.

\renewcommand{\labelenumi}{(\roman{enumi})}
\begin{enumerate}

\item
There exists
$1 \leq j \leq s$ such that $\deg \shR_j < 0$
and $c_j \neq 0$, where $c_j $ denotes the component of $c$
in $H^1(Y, \shR_j)$.

\item
The complement of the forcing divisor is affine.

\item
$f_0 \not\in (f_1, \ldots ,f_n)^* $.

\end{enumerate}
\end{corollary}
\proof
The equivalence (ii) $ \Leftrightarrow $ (iii) was stated in
\ref{geokrittightplus},
and (i) $\Leftrightarrow$ (ii) is \ref{numkritaffin}.
\qed

\begin{corollary}
\label{numkritplus}
Let $K$ be an algebraically closed field of positive characteristic
and let
$R$ be a normal homogeneous coordinate ring over the
elliptic curve $Y= \Proj \, R$.
Let $f_1, \ldots ,f_n$ be homogeneous generators of
an $R_+$-primary homogeneous ideal in $R$
and let $m \in \NN$ be a number.
Let $\shR(m)$ be the corresponding locally free sheaf on $Y$
and let $\shR(m) = \shR_1 \oplus \ldots \oplus \shR_s$
be the decomposition into indecomposable locally free sheaves.
Let $f_0 \in R$ be another homogeneous element of degree $m$
and let $c \in H^1(Y, \shR(m))$ be the corresponding class.
Then the following are equivalent.

\renewcommand{\labelenumi}{(\roman{enumi})}
\begin{enumerate}

\item
There exists
$1 \leq j \leq s$ such that $\deg \shR_j < 0$
and $c_j \neq 0$, where $c_j $ denotes the component of $c$
in $H^1(Y, \shR_j)$.

\item
The forcing divisor intersects every curve.

\item
The complement of the forcing divisor is affine.

\item
$f_0 \not\in (f_1, \ldots ,f_n)^* $.

\item
$f_0 \not\in (f_1, \ldots , f_n)^{gr +}$.

\end{enumerate}
\end{corollary}
\proof
This follows directly from \ref{numkritcurve}, \ref{numkrittight}
and \ref{geokrittightplus}.
\qed

\medskip
Our main theorem is now easy to deduce.

\begin{theorem}
\label{tightplus}
Let $K$ be an algebraically closed field of positive characteristic
and let
$R$ denote a normal homogeneous coordinate ring over an elliptic curve.
Let $\foa $ be an $R_+$-primary homogeneous ideal in $R$.
Then $\foa ^{gr +} = \foa^+ = \foa^* $.

This holds in particular for $R=K[x,y,z]/(F)$,
where $F$ is homogeneous of degree $3$ and $R$ is normal.
\end{theorem}
\proof
The inclusions $ \subseteq $ are clear. It is known that the tight closure
of a homogeneous ideal is again homogeneous,
see \cite[Theorem 4.2]{hochsterhunekesplitting}.
Hence the statement follows from \ref{numkritplus}.
\qed

\begin{remark}
If the $p$-rank (=Hasse invariant) of the elliptic curve $Y$ is 0,
then the plus closure (=tight closure) of a primary homogeneous ideal is
the same as its Frobenius closure.
This follows from the proof of Lemma \ref{pvanish} and the proof of
(iii) $\Rightarrow$ (i) of Theorem \ref{numkritcurve}.
\end{remark}

\medskip
The tight closure of a primary homogeneous ideal is
easy to describe by a numerical condition,
if (there exists a system of homogeneous generators such that)
the corresponding
vector bundle is indecomposable.

\begin{corollary}
\label{indecomposable}
Let $K$ be an algebraically closed field and let
$R$ be a normal homogeneous coordinate ring over the
elliptic curve $Y= \Proj \, R$.
Let $f_1, \ldots ,f_n$ be homogeneous generators of
an $R_+$-primary homogeneous ideal in $R$
with $\deg \, f_i = d_i$ and let $m \in \NN$ be a number.
Suppose that the corresponding locally free sheaf $\shR(m)$ on $Y$
is indecomposable.
Let $f_0 \in R$ be a homogeneous element of degree $m$,
defining the cohomology class $c \in H^1(Y,\shR(m))$.
Then the following are equivalent.

\renewcommand{\labelenumi}{(\roman{enumi})}
\begin{enumerate}

\item
$d_1+ \ldots +d_n -(n-1)m >0$ and $f_0 \not \in (f_1, \ldots ,f_n)$.

\item
$\deg \shR < 0$ and the cohomology class is $c \neq 0$.

\item
The sheaf $\shF(-m)=\shR(m)^\dual$ is ample and $c \neq 0$.

\item
The sheaf $\shF'(-m)$ is ample.

\item
The complement of the forcing divisor is affine.

\item
$f_0 \not\in (f_1, \ldots ,f_n)^* \,
{\rm(}= (f_1, \ldots ,f_n)^+ $ in positive
characteristic{\rm )}.

\end{enumerate}

\end{corollary}

\proof
The degree of $\shR(m)$ is $\deg \shR(m) =-3(d_1+ \ldots +d_n -(n-1)m)$
due to \ref{buendelalsproj}, hence (i) $\Leftrightarrow$ (ii) is clear.
(ii) $\Leftrightarrow$ (iii) follows from \ref{amplekrit},
(iii) $\Rightarrow $ (iv) follows from \cite[Proposition 2.2]{giesekerample}
and \ref{amplekrit}, (iv) $\Rightarrow $ (v) is clear
and (v) $\Rightarrow $ (ii) is \ref{numkrittight} for indecomposable $\shR(m)$.
\qed

\begin{corollary}
\label{indecomposable2}
Let $K$ be an algebraically closed field and let
$R$ be a normal homogeneous coordinate ring over the
elliptic curve $Y= \Proj \, R$.
Let $f_1, \ldots ,f_n$ be homogeneous generators of
an $R_+$-primary homogeneous ideal in $R$
with $\deg \, f_i = d_i$ and suppose
that the corresponding locally free sheaf $\shR(m)$ on $Y$
is indecomposable.
Set $k= \frac{d_1 + \ldots +d_n}{n-1} $.
Then
$$(f_1, \ldots ,f_n)^* = (f_1, \ldots ,f_n) + R_{\geq k} $$
\end{corollary}

\proof
Let $f_0 \in R$ be homogeneous of degree $m$.
Suppose first that $m \geq k$.
Then $(n-1)m \geq d_1 + \ldots +d_n$ and hence
the numerical condition in \ref{indecomposable}
is not fulfilled, thus $f_0 \in (f_1, \ldots ,f_n)^*$.

Suppose now that $f \in (f_1, \ldots ,f_n)^*$.
Then from \ref{indecomposable} we see that
$d_1+ \ldots +d_n -(n-1)m \leq 0$ or $f_0 \in (f_1, \ldots ,f_n)$,
which gives the result.
\qed

\begin{example}
Let the elliptic curve $Y$ be given by the equation $x^3+y^3+z^3=0$.
Then the ideal $(x^2,y^2,z^2)$ defines an indecomposable sheaf of
syzygies. Look at the total degree $m=3$. Then $\shR(3) $ has degree
$0$ and its determinant is trivial. The syzygy $(x,y,z)$ gives a
global section (up to multiples the only section) of $\shR(3)$. The
cokernel of $0 \ra \O_Y \ra \shR(3)$ is again invertible (check
locally), hence isomorphic to $\O_Y$. Furthermore the sequence does
not split, hence $\shR(3)$ is indecomposable. In particular it
follows that $xyz \in (x^2,y^2,z^2)^*$.
\end{example}

\begin{example}
The statements in \ref{indecomposable} and \ref{indecomposable2} do
not hold without the condition that $\shR(m)$ is indecomposable. The
easiest way to obtain decomposable sheaves of syzygies is to look at
redundant systems of generators. Consider again the elliptic curve
$Y$ given by $x^3+y^3+z^3=0$.

Look at the elements $x^2,y^2,x^2$.
Then the corresponding sheaf is of course decomposable.
For $m=3$, the sheaf
$\shR(3) \subset \O_Y(1) \oplus \O_Y(1) \oplus \O_Y(1) $
is given by $(g_1,g_2,g_3)$ such that $ g_1x^2+g_2y^2+g_3x^2=0$
and it is easy to see that
$\shR(3) \cong \O_Y(1) \oplus \O_Y(-1)$.
Let $f_0 =xyz$. Then the number $d_1+ \ldots +d_n -(n-1)m $ is zero, but
$xyz \not\in (x^2,y^2,x^2)^\pasoclo =(x^2,y^2)^\pasoclo$.
The complement of the forcing divisor in $\PP(x^2,y^2,x^2,xyz)$
is affine, but it is not ample, since its degree is zero (or
since $\O_Y(-1)$ is a quotient invertible sheaf
of $\shF'(-3)$ of negative degree).

Consider now $x,y,z^3$ and $f_0=z^2$. Then $z^2 \in
(x,y,z^3)^*=(x,y)^*$, but $z^2 \not\in (x,y)$ and the number in
\ref{indecomposable}(i) is $1+1+3-2\cdot 2=1 >0$. The sheaf of
syzygies decomposes $\shR(3)= \O_Y \oplus \O_Y(1)$.

The ideal $(x^2,y^2,xy)$ provides an example where no generator is
superfluous, but the corresponding sheaf of syzygies is anyway
decomposable.
\end{example}

\medskip
Recall that for a locally free sheaf the slope is defined by
$\mu (\shF)= \deg\, \shF / \rank \shF$. Furthermore
$ \mu_{\rm min} (\shF)= {\rm min} \{ \mu (\shG): \, \shF \ra \shG \ra 0 \}$
and
$ \mu_{\rm max} (\shF)= {\rm max} \{ \mu (\shS): \, 0 \ra \shS \ra \shF \}$.
In the case of an elliptic curve it is easy to see that
$\mu_{\rm min} (\shF)= {\rm min}_j \, \mu (\shF_j)$ and
$\mu_{\rm max} (\shF)= {\rm max}_j \, \mu (\shF_j)$, where
$\shF =\shF_1 \oplus \ldots \oplus \shF_s$ is the decomposition
in indecomposable sheaves.

\begin{corollary}
\label{maxmin} Let $K$ be an algebraically closed field and let
$R=K[x,y,z]/(F)$, where $F$ is a homogeneous polynomial of degree
$3$ defining the elliptic curve $Y= \Proj \, R$. Let $f_1, \ldots
,f_n$ be homogeneous generators of a primary graded ideal in $R$,
$d_i= \deg \, f_i$. Let $\shR(0)$ be the corresponding locally free
sheaf of syzygies on $Y$ of total degree $0$ and let $\shF(0)$ be
its dual sheaf. Let $f_0$ denote another homogeneous element of
degree $m$. Then the following hold.

\renewcommand{\labelenumi}{(\roman{enumi})}
\begin{enumerate}

\item
If $m \geq  \frac{1}{3}  \mu_{\rm max} (\shF(0))$,
then $f_0 \in (f_1, \ldots,f_n)^*$.

\item
If $m < \frac{1}{3} \mu_{\rm min} (\shF(0))$,
then $f_0 \in (f_1, \ldots, f_n)^*$
if and only if $f_0 \in (f_1, \ldots ,f_n)$.

\item
If $\shF(0)$ is semistable, then
$(f_1, \ldots ,f_n)^* = (f_1, \ldots ,f_n) + R_{\geq k} $,
where $k= \frac{d_1 + \ldots +d_n}{n-1}$.

\end{enumerate}
\end{corollary}
\proof
Let $\shR(0)=\shR_1 \oplus \ldots \oplus \shR_s$ be the decomposition into
indecomposable locally free sheaves.
The homogeneous element $f_0$ of degree $m$ defines a cohomology class
$c \in H^1(Y,\shR \otimes \O_Y(m))$.

(i).
The condition is that $3m \geq \mu(\shF_j)=-\mu (\shR_j)$ holds
for every $j$. This means that
$$\deg\, (\shR_j \otimes \O_Y(m))
= \deg\, (\shR_j) + 3m \rank \,(\shR_j)
\geq \deg\, (\shR_j) - \mu (\shR_j) \rank \,( \shR_j) = 0 \, .$$
Hence the numerical condition in \ref{numkrittight}
is not fulfilled and $f_0 \in (f_1, \ldots ,f_n)^*$.

(ii).
In this case we have $ 3m < \mu(\shF_j) = - \mu (\shR_j)$ for all $j$, hence
$\deg \, (\shR_j \otimes \O_Y(m)) < 0$,
and the numerical criterion in \ref{numkrittight}
is true if and only if the cohomological class is $c \neq 0$.
This gives the result.

(iii).
The sheaf $\shF(0)$ is semistable if and only if
$\mu_{\rm max} (\shF(0))= \mu_{\rm min} (\shF(0))$.
This number equals then the slope
$\mu (\shF(0)) = 3 (d_1 + \ldots +d_n)/(n-1)$,
so the result follows from (i) and (ii).
\qed

%===========================================================

\end{document}